\documentclass[12pt,a4paper]{article}
\usepackage{pifont}

\usepackage{pst-node}

\usepackage{epsfig}
\usepackage{amsmath}
\usepackage{amssymb}
\usepackage{graphicx}

\usepackage{tikz}

\newtheorem{propo}{Proposition}[section]
\newtheorem{defi}[propo]{Definition}

\newtheorem{lemma}[propo]{Lemma}

\newtheorem{theo}[propo]{Theorem}

\newtheorem{rem}[propo]{Remark}

\newcommand{\bl}{\begin{lemma}}
\newcommand{\el}{\end{lemma}}

\def\d12{{_{12}}}

\def\Cay{{\rm Cay}}

\def\Cos{{\rm Cos}}

\usepackage{indentfirst,latexsym,bm}

\def\PSL{{\rm PSL}}

\def\Aut{{\rm Aut}}

\def\diam{{\rm diam}}

\def\girth{{\rm girth}}

\begin{document}
\title{Finite 2-geodesic transitive graphs of  prime valency
}

\author{Alice Devillers$^{1}$, Wei Jin$^{2}$\thanks{Corresponding author is supported by  NNSF (11301230), Jiangxi (GJJ14351, 20142BAB211008) and UWA (SIRF).}, Cai
Heng Li$^{1}$ and Cheryl E. Praeger$^{1}$\thanks{ This paper  forms part of
Australian Research Council Federation Fellowship FF0776186 held by
the fourth author. The fourth author is also affiliated with King
Abdulazziz University, Jeddah. The first author is supported by UWA
as part of the Federation Fellowship project. }\thanks{E-mail
addresses: alice.devillers@uwa.edu.au(A. Devillers),
jinweipei82@163.com(W. Jin), cai.heng.li@uwa.edu.au(C. H.
Li) and cheryl.praeger@uwa.edu.au(C. E. Praeger).  }
\\{\footnotesize
$^{1}$Centre for the Mathematics of Symmetry and Computation, School of
Mathematics and Statistics,  }
\\{\footnotesize
The University of Western Australia, Crawley, WA 6009, Australia}
\\{\footnotesize
$^{2}$School of Statistics, Research Centre of Applied Statistics  }
\\{\footnotesize
Jiangxi University of Finance and Economics,  Nanchang, Jiangxi, 330013, P.R.China  }}

\date{ }

\maketitle

\begin{abstract}

We classify non-complete prime valency  graphs satisfying the
property  that their automorphism group is transitive on both the
set of arcs and the set of  $2$-geodesics. We prove that either
$\Gamma$ is  2-arc transitive or  the valency $p$ satisfies $p\equiv
1\pmod 4$, and for each such prime there is a unique graph with this
property: it is a non-bipartite antipodal double cover of the
complete graph $K_{p+1}$ with automorphism group $PSL(2,p)\times
Z_2$ and diameter 3.

\end{abstract}

\vspace{2mm}

 \hspace{-17pt}{\bf Keywords:} 2-geodesic transitive graph; 2-arc transitive graph; cover

\section{Introduction}

In this paper, graphs are finite, simple and undirected. For a graph
$\Gamma$,   a vertex triple $(u,v,w)$  with $v$ adjacent to both $u$
and $w$ is called a \emph{$2$-arc} if $u\neq w$, and  a
\emph{$2$-geodesic} if in addition $u,w$ are not adjacent. An
\emph{arc} is an ordered pair of adjacent vertices. A non-complete
graph $\Gamma$   is said to be \emph{ $2$-arc transitive} or \emph{
$2$-geodesic transitive} if its automorphism group  is transitive on
arcs, and also on  2-arcs or 2-geodesics, respectively. Clearly, every
2-geodesic is a 2-arc, but some 2-arcs may not be 2-geodesics. If
$\Gamma$ has girth 3 (length of the shortest cycle is 3), then the
2-arcs contained in $3$-cycles are not 2-geodesics. Thus the family
of non-complete $2$-arc transitive graphs is properly contained in
the family of $2$-geodesic transitive graphs. The graph in Figure 1
is the icosahedron which is 2-geodesic transitive but not 2-arc
transitive with valency 5.

\medskip


\begin{figure}\label{fig1}
\centering

\begin{tikzpicture}

\draw (1,3)-- (3,3); \draw (3,3)-- (5,3); \draw (5,3)-- (7,3);

\draw (1,3)-- (3,1); \draw (3,1)-- (5,1); \draw (5,1)-- (7,3);

\draw (1,3)-- (3,2); \draw (3,2)-- (5,2); \draw (5,2)-- (7,3);

\draw (1,3)-- (3,4); \draw (3,4)-- (5,4); \draw (5,4)-- (7,3);

\draw (1,3)-- (3,5); \draw (3,5)-- (5,5); \draw (5,5)-- (7,3);

\draw (3,1)-- (3,2); \draw (3,2)-- (3,3); \draw (3,3)-- (3,4); \draw
(3,4)-- (3,5);

\draw (5,1)-- (5,2); \draw (5,2)-- (5,3); \draw (5,3)-- (5,4); \draw
(5,4)-- (5,5);

\draw (3,5)-- (5,5); \draw (3,5)-- (5,4);

\draw (3,4)-- (5,4); \draw (3,4)-- (5,3);

\draw (3,3)-- (5,3); \draw (3,3)-- (5,2);

\draw (3,2)-- (5,2); \draw (3,2)-- (5,1);

\draw (3,1)-- (5,1); \draw (3,1)-- (5,5);

\draw (3,1) .. controls (3.5,1) and (3.5,5) .. (3,5);

\draw (5,1) .. controls (4.5,1) and (4.5,5) .. (5,5);


\filldraw[black] (1,3) circle (2pt) (3,3) circle (2pt) (5,3) circle
(2pt) (7,3) circle (2pt);

\filldraw[black] (3,1) circle (2pt) (5,1) circle (2pt);

\filldraw[black] (3,2) circle (2pt) (5,2) circle (2pt);

\filldraw[black] (3,4) circle (2pt) (5,4) circle (2pt);

\filldraw[black] (3,5) circle (2pt) (5,5) circle (2pt);

\end{tikzpicture}
\caption{Icosahedron}
\end{figure}
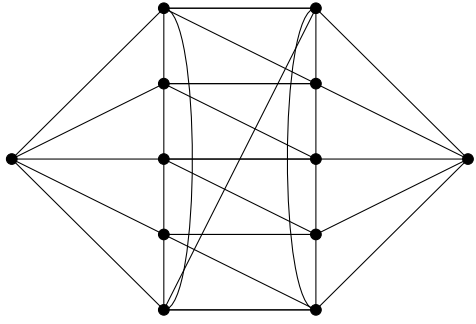

The study of  2-arc transitive graphs goes back to  Tutte
\cite{Tutte-1,Tutte-2}. Since then, this family of graphs has been
studied extensively, see
\cite{Baddeley-1,IP-1,Praeger-4,Weiss-1,weiss}. In this paper, we
are interested in $2$-geodesic transitive graphs, in particular,
which are not 2-arc transitive, that is, they have girth 3.  We
first construct a family of coset graphs, and  prove that each of
these graphs   is  2-geodesic transitive but not 2-arc transitive of
prime valency. We then prove that each graph with these properties
belongs to the family.

For a finite group $G$,  a  core-free  subgroup $H$ (that is,
$\bigcap_{g\in G}H^g=1$), and an element $g\in G$ such that
$G=\langle H,g\rangle$ and $HgH=Hg^{-1}H$, the \emph{coset graph}
$\Cos(G,H,HgH)$ is the  graph with vertex set $\{Hx|x\in G\}$, such
that  two vertices $Hx,Hy$  are adjacent if and only if $yx^{-1}\in
HgH$. This graph is   connected, undirected, and $G$-arc transitive
of valency $|H:H\cap H^g|$, see \cite{Lorimer-1}.

\begin{defi}\label{coset-constr}
{\rm Let $\mathcal{C}(5)$ be the singleton set containing the
icosahedron, and for a prime $p>5$ with $p\equiv 1\pmod{4}$, let
$\mathcal{C}(p)$  consist of the  coset graphs $\Cos(G,H,HgH)$ as
follows.   Let $G=PSL(2,p)$, choose $a\in G$ of order $p$, so
$N_{G}(\langle a\rangle)=\langle a\rangle:\langle b\rangle\cong
Z_p:Z_{\frac{p-1}{2}}$ for some $b\in G$ of order $\frac{p-1}{2}$.
Then $N_{G}(\langle b^2\rangle)=\langle b\rangle:\langle
c\rangle\cong D_{p-1}$ for some $c\in G$ of order 2. Let $H=\langle
a\rangle:\langle b^2\rangle$ and $g=cb^{2i}$ for some $i$.  }
\end{defi}

These graphs have appeared a number of times in the literature. They
were constructed by D. Taylor \cite{Taylor-1} as a family of regular
two-graphs (see \cite[p.14]{BCN}), they  appeared in the
classification of antipodal distance transitive covers of complete
graphs in \cite{GLP}, and were also    constructed explicitly as
coset graphs and studied by the third author in
\cite{Li-ci-soluble}. (Antipodal covers of graphs are defined in
Section 2.)

A path of shortest length from a vertex $u$ to a vertex $v$ is
called a \emph{geodesic} from $u$ to $v$, or sometimes an
\emph{$i$-geodesic} if the distance between $u$ and $v$ is $i$.  The
graph $\Gamma$ is said to be \emph{geodesic transitive} if  its
automorphism group is transitive on the set of $i$-geodesics for all
positive integers $i$ less than or equal to the diameter of
$\Gamma$.

\begin{theo}\label{gt-primeval2}
{\rm (a)} A graph    $\Gamma \in \mathcal{C}(p)$ if and only if
$\Gamma$ is a connected non-bipartite antipodal double cover of
$K_{p+1}$ with $p\equiv 1 \pmod 4$, and $\Aut \Gamma \cong
PSL(2,p)\times Z_2$.

{\rm (b)} For a given $p$, all graphs in $\mathcal{C}(p)$ are
isomorphic, geodesic transitive and  have diameter $3$.
\end{theo}

Our second  result shows that the graphs in Definition
\ref{coset-constr} are the only 2-geodesic transitive graphs of
prime valency that are not 2-arc transitive.

\begin{theo}\label{gt-primeval3}
Let $\Gamma$ be a connected non-complete  graph    of prime valency
$p$. Then $\Gamma$ is $2$-geodesic transitive if and only if
$\Gamma$ is $2$-arc transitive, or $p\equiv  1\pmod 4$ and $\Gamma
\in \mathcal{C}(p)$.
\end{theo}

These  two theorems show that up to isomorphism, there is a unique
connected 2-geodesic transitive but not 2-arc transitive graph of
prime valency $p$ and $p\equiv 1 \pmod 4$. The family of 2-geodesic
transitive but not 2-arc transitive graphs of valency 4 has been
determined in \cite{DJLP-2}.   It would be interesting to know if a
similar classification is possible for non-prime valencies at least
6.  This is the subject of further research by the second author, see \cite{JLX-val6}.

\section{Preliminaries}

In this section, we give some definitions and prove some results
which will be used in the following discussion.  Let $\Gamma$ be a
graph. We use $V\Gamma$,  $E\Gamma$ and $\Aut \Gamma$ to denote its
\emph{vertex set},  \emph{edge set} and \emph{automorphism group},
respectively. The size of $V\Gamma$ is called the \emph{order} of
the graph.  The graph $\Gamma$ is said to be \emph{vertex
transitive}  if the action of $\Aut \Gamma$ on $V\Gamma$  is
transitive.

For two distinct vertices $u,v$ of $\Gamma$, the smallest value for
$n$ such that there is a path of length $n$ from $u$ to $v$ is
called the \emph{distance} from $u$ to $v$ and is denoted by
$d_{\Gamma}(u, v)$.  The \emph{diameter} $\diam(\Gamma)$ of a
connected graph $\Gamma$ is the maximum of $d_{\Gamma}(u, v)$ over
all $u, v \in V\Gamma$. We set $\Gamma_2(v)=\{u\in
V\Gamma|d_{\Gamma}(v,u)=2\}$ for every vertex $v$.

\medskip

Quotient graphs play an important role in this paper. Let $G$ be a
group of permutations acting on a set $\Omega$. A
\emph{$G$-invariant partition} of $\Omega$  is a partition
$\mathcal{B}=\{B_1,B_2,\ldots,B_n\}$ such that for each $g\in G$,
and each $B_i \in \mathcal{B}$, the image $B_{i}^g\in \mathcal{B}$.
The parts of $\Omega$ are often called \emph{blocks} of $G$ on
$\Omega$. For a $G$-invariant partition $\mathcal{B}$ of $\Omega$,
we have two smaller transitive permutation groups, namely the  group
$G^{\mathcal{B}}$ of permutations of $\mathcal{B}$ induced by $G$;
and the group $G_{B_i}^{B_i}$ induced on $B_i$ by $G_{B_i}$ (the setwise stabiliser of $B_i$ in $G$) where
$B_i\in \mathcal{B}$.  Let $\Gamma$ be a graph, and let $G\leq \Aut
\Gamma$. Suppose $\mathcal{B}=\{B_1,B_2,\ldots,B_n \}$ is a
$G$-invariant partition of $V\Gamma$. The \emph{quotient graph}
$\Gamma_{\mathcal{B}}$ of $\Gamma$ relative to $\mathcal{B}$ is
defined to be the graph with vertex set $\mathcal{B}$ such that
$\{B_i,B_j\}$ ($i\neq j$) is an edge of $\Gamma_{\mathcal{B}}$ if and only if
there exist $x\in B_i, y\in B_j$ such that $\{x,y\} \in E\Gamma$. We
say that $\Gamma_{\mathcal{B}}$ is \emph{nontrivial} if $1<
|\mathcal{B}|< |V\Gamma|$.  The graph $\Gamma$ is said to be a
\emph{cover} of $\Gamma_{\mathcal{B}}$ if for each edge
$\{B_i,B_j\}$ of $\Gamma_{\mathcal{B}}$ and $v\in B_i$, we have
$|\Gamma(v)\cap B_j|=1$.

\medskip

For a graph $\Gamma$, the \emph{k-distance graph} $\Gamma_k$ of
$\Gamma$ is the graph with vertex set $V\Gamma$, such that  two
vertices are adjacent if and only if they are at distance $k$ in
$\Gamma$. If $d=\diam(\Gamma)\geq 2$ and $\Gamma_d$ is a disjoint
union of complete graphs, then  $\Gamma$ is said to be an
\emph{antipodal graph}. In other words, the vertex set of an
antipodal graph $\Gamma$ of diameter $d$, may be partitioned into
so-called  \emph{fibres}, such that any two distinct vertices in the
same fibre are at distance $d$ and two vertices in different fibres
are at distance less than $d$.      For an antipodal  graph $\Gamma$
of diameter $d$, its \emph{antipodal quotient graph} $\Sigma$ is the
quotient graph of $\Gamma$ where $\mathcal{B}$ is the set of fibres.
If further, $\Gamma$ is a cover of $\Sigma$, then $\Gamma$ is called
an \emph{antipodal cover} of $\Sigma$.

\medskip

Paley graphs were first defined by Paley in 1933, see
\cite{Paley-1}. These graphs are vertex transitive, self-complementary,  and have many
nice properties. Let $q=p^e$ be a prime power such that $q\equiv 1
\pmod{4}$. Let $F_q$ be the finite field of order $q$. The
\emph{Paley graph} $P(q)$ is the graph with vertex set $F_q$, where
two distinct vertices $u,v$ are adjacent if and only if $u-v$ is a
nonzero square in $F_q$. The congruence condition on $q$ implies
that $-1$ is a square in $F_q$, and hence $P(q)$ is an undirected
graph.


Lemma \ref{arctrans-order-p-1}  is used in the proof of Theorem
\ref{gt-primeval3}, and its proof uses the following famous result
of Burnside.

\begin{lemma}{\rm( \cite[Theorem 3.5B]{DM-1})}\label{val-p-prim-1}
A  primitive permutation group $G$ of prime degree $p$  is either
$2$-transitive, or solvable and $G\leq AGL(1,p)$.
\end{lemma}

For a finite group $G$, and a subset $S$ of $G$ such that $1\notin
S$ and $S=S^{-1}$, the \emph{Cayley graph} $\Cay(G,S)$ of $G$ with
respect to $S$ is  the graph with vertex set $G$ and edge set
$\{\{g,sg\} \,|\,g\in G,s\in S\}$. The Paley graph $P(q)$ is a
Cayley graph for the additive group $G=F_{q}^+$ with
$S=\{w^2,w^4,\ldots,w^{q-1}=1\}$, where  $w$ is a primitive element
of $F_q$.

\begin{lemma}\label{arctrans-order-p-1}
Let $\Gamma$ be an arc transitive graph of  prime order $p$ and
valency $\frac{p-1}{2}$. Then $p\equiv 1 \pmod 4$, $\Aut \Gamma
\cong Z_p:Z_{\frac{p-1}{2}}$, and  $\Gamma \cong P(p)$.
\end{lemma}
{\bf Proof.} Since $\Gamma$ has valency $\frac{p-1}{2}$, $p$ is an
odd prime. Since $\Gamma$ has the given order and valency, it
follows that $\Gamma$ has $p(\frac{p-1}{2})/2$ edges. This implies
that  $p\equiv 1 \pmod 4$.

Let $A=\Aut \Gamma$. Since  $A$ is transitive on $V\Gamma$ and  $p$
is a prime, $A$ is primitive on $V\Gamma$, and since $\Gamma$ is arc
transitive, $|A|$ is divisible by $\frac{p(p-1)}{2}$. Since $\Gamma$
is neither complete nor empty, it follows by Lemma
\ref{val-p-prim-1} that $A< AGL(1,p)=Z_p:Z_{p-1}$. Thus $|A|$ is a
proper divisor of $p(p-1)$, and at least $\frac{p(p-1)}{2}$, and so
$|A|=\frac{p(p-1)}{2}$. Hence   $A \cong Z_p:Z_{\frac{p-1}{2}}$.

Since $Z_p$ is regular on $V\Gamma$, it follows from \cite[Lemma
16.3]{Biggs-1} that $\Gamma$ is a Cayley graph for $Z_p$. Thus
$\Gamma=\Cay(G,S)$ where $G\cong Z_p$, $S\subseteq G\setminus
\{0\}$, $S=S^{-1}$ and $|S|=\frac{p-1}{2}$. Now we may identify $G$
with $F_p^+$ where $F_p$ is a finite field of order  $p$. Let $v\in
V\Gamma$ be the vertex corresponding to $0\in G$.  Then $A_v$ is the
unique subgroup of order $\frac{p-1}{2}$ of $F_p^*=\langle
w\rangle$, that is, $A_v=\langle w^2\rangle$. The $A_v$-orbits in
$F_p$ are $\{0\}$, $S_1=\{w^2,w^4,\ldots,w^{p-1}\}$ and
$S_2=\{w,w^3,\ldots,w^{p-2}\}$, and so $S=S_1$ or $S_2$, and
$\Gamma=P(p)$ or its complement respectively. In either case,
$\Gamma\cong P(p)$.  $\Box$

\medskip

To end the section, we cite a property of Paley graphs which will be
used in the next section.

\begin{lemma}{\rm(\cite[p.221]{GR})}\label{val-p-rem1}
Let  $\Gamma=P(q)$, where  $q$ is a prime power such that $q\equiv 1
\pmod 4$. Let $u,v$ be   distinct vertices of $\Gamma$. If $u,v$ are
adjacent, then $|\Gamma(u)\cap \Gamma(v)|=\frac{q-5}{4}$; if $u,v$
are not adjacent, then $|\Gamma(u)\cap \Gamma(v)|=\frac{q-1}{4}$.
\end{lemma}

\section{Proof of Theorem \ref{gt-primeval2} }

We study graphs in the family $\mathcal{C}(p)$ for each prime
$p\equiv 1 \pmod 4$. We first  collect some properties of graphs in
$\mathcal{C}(p)$ for $p>5$, which can be found in \cite[Theorem
1.1]{Li-ci-soluble} and its proof.

\begin{rem}\label{smallval-con-rem}
{\rm Let  $\Gamma \in \mathcal{C}(p)$ and $p>5$. Then $G=\langle
H,g\rangle$, $\Gamma$ is connected and $G$-arc transitive of valency
$p$, $\Aut \Gamma \cong G\times Z_2$, $|V\Gamma|=|G:H|=2p+2$.
Further, $\diam(\Gamma)=\girth(\Gamma)=3$, so    $\Gamma$ is  not
$2$-arc transitive.

The orbit set $\mathcal{B}=\{\Delta_1,\Delta_2,\ldots,\Delta_{p+1}\}$ of the normal subgroup $K\cong Z_2$ of $\Aut \Gamma$ forms a system of imprimitivity for $\Aut \Gamma$ in $V \Gamma$, and it follows from the proof of \cite[Theorem 1.1]{Li-ci-soluble} that this is the unique nontrivial system of imprimitivity and the kernel of the action of
$\Aut \Gamma$ on $\mathcal{B}$ is the normal subgroup $K$. For $i=1,\ldots,p+1$, let $\Delta_i=\{v_i,v_i'\}$. Then  $v_i$ is not
adjacent to $v_i'$, and for each $j\neq i$, $v_i$ is adjacent to
exactly one point of $\Delta_j$ and $v_i'$ is adjacent to the other.
Thus, $\Gamma(v_1)\cap \Gamma(v_1')=\emptyset$, $V\Gamma=\{v_1\}\cup
\Gamma(v_1) \cup \{v_1'\}\cup \Gamma(v_1') $, and $\Gamma$ is a
non-bipartite double cover of $K_{p+1}$.
}
\end{rem}

The next lemma  shows that graphs in $\mathcal{C}(p)$ are geodesic
transitive.

\begin{lemma}\label{smallval-p}
Let $p$ be a prime and $p\equiv 1 \pmod 4$. Then each graph in
$\mathcal{C}(p)$  is geodesic transitive of girth $3$ and diameter
$3$.
\end{lemma}
{\bf Proof.} Let $\Gamma\in \mathcal{C}(p)$. If $p=5$, then $\Gamma$
is the icosahedron of girth 3 and diameter 3. Its automorphism group
is $PSL(2,5)\times Z_2$ and it is geodesic transitive. Now suppose
that $p>5$.  Let $\mathcal {B}$ be as in Remark
\ref{smallval-con-rem}, $A:=\Aut \Gamma$, $v_1\in V\Gamma$ and $u\in
\Gamma(v_1)$. Let $K$ be the kernel of the $A$-action on $\mathcal
{B}$ so that the induced group  $A^{\mathcal {B}}=A/K$. Then by the
proof of \cite[Theorem 1.1]{Li-ci-soluble}, $K\cong Z_2 \lhd A$,
$A=G\times K$, $A^{\mathcal {B}}\cong G=PSL(2,p)$ and $(A^{\mathcal
{B}})_{\Delta_1}\cong A_{v_1}$.  Since $A\cong G\times Z_2$, it
follows that $|A_{v_1}|=\frac{p(p-1)}{2}$, and by Lemma 2.4 of
\cite{Li-ci-soluble}, $A_{v_1}\cong Z_p:Z_{\frac{p-1}{2}}$, which
has a unique permutation action of degree $p$, up to permutational
isomorphism. Since $\Gamma$ is $A$-arc transitive, $A_{v_1}$ is
transitive on $\Gamma(v_1)$ and hence on $\mathcal{B}\setminus
\{\Delta_1\}$, and therefore also on $\Gamma(v_1')$, all of degree
$p$. Thus the $A_{v_1}$-orbits in $V\Gamma$ are
$\{v_1\},\Gamma(v_1),\Gamma(v_1')$ and $\{v_1'\}$, and it follows
that $\Gamma(v_1')=\Gamma_2(v_1)$. Moreover, $A_{v_1,u}\cong
Z_{\frac{p-1}{2}}$ has orbit lengths $1,\frac{p-1}{2},\frac{p-1}{2}$
in $\Gamma(v_1)$, and hence has the same orbit lengths in
$\Gamma_2(v_1)$, and also in $\Gamma(u)$ (since $A_{v_1,u}$ is the
point stabiliser of $A_u$ acting on $\Gamma(u)$). Since
$\Gamma(v_1)\cap \Gamma(u)\neq \emptyset$, it follows that the
$A_{v_1,u}$-orbits in $\Gamma(u)$ are $\{v_1\},\Gamma(v_1)\cap
\Gamma(u)$, and $\Gamma_2(v_1)\cap \Gamma(u)$. Thus  $\Gamma$ is
$(A,2)$-geodesic transitive and $\girth(\Gamma)=3$. Further, as
$\Gamma_3(v_1)=\{v_1'\}$, it follows that  $\Gamma$ is geodesic
transitive and has diameter 3.  $\Box$

\medskip

In the proof of the second part of  Theorem \ref{gt-primeval2}, we
repeatedly use the fact that each $\sigma \in \Aut G$ induces an
isomorphism from $\Cos(G,H,HgH)$ to
$\Cos(G,H^{\sigma},H^{\sigma}g^{\sigma}H^{\sigma})$, and in
particular, we use this fact for the conjugation action by elements
of $G$. For a subset $\Delta$ of the vertex set of a graph $\Gamma$,
we use  $[\Delta]$ to denote  the subgraph of $\Gamma$ induced by
$\Delta$.

\bigskip

\noindent{\bf Proof of Theorem \ref{gt-primeval2}} (a) Suppose first
that $\Gamma$ is a connected non-bipartite antipodal double cover of
$K_{p+1}$ with $p\equiv 1 \pmod 4$, and $A:=\Aut \Gamma \cong
PSL(2,p)\times Z_2$. Then  $|V\Gamma|=2p+2$, and for each  $u\in
V\Gamma$,  let  $u'\in V\Gamma$ be its unique vertex at maximum
distance. Then  $|\Gamma(u)|=p=|\Gamma(u')|$, and $\Gamma(u)\cap
\Gamma(u')=\emptyset$. Since $\Gamma$ is connected, it follows that
$V\Gamma=\{u\}\cup \Gamma(u)\cup \Gamma(u')\cup \{u'\}$, and  the
diameter of $\Gamma$ is 3.

Let $\mathcal{B}=\{B_1,B_2,\ldots,B_{p+1}\}$ be the invariant
partition of $V\Gamma$ such that $\Gamma_{\mathcal{B}}\cong K_{p+1}$
and $\Gamma$ is a  non-bipartite antipodal double cover of
$\Gamma_{\mathcal{B}}$. Let $K$ be the kernel of the $A$-action on
$\mathcal{B}$. As each $|B_i|=2$, it follows that $K$ is a 2-group.
Further, as $K$ is a normal subgroup of $A$ and $PSL(2,p)$ is a
simple group, it follows that $K\cong Z_2$. Thus $G:=PSL(2,p)$ acts
faithfully on $\mathcal{B}$. Since the $G$-action on $p+1$ points is
unique and this action is 2-transitive, it follows that $G$ is
2-transitive on $\mathcal{B}$, and so $\Gamma_{\mathcal{B}}$ is
$G$-arc transitive. Thus either $G$ is transitive on $V\Gamma$ or
$G$ has two orbits $\Delta_1,\Delta_2$ in $V\Gamma$ of size $p+1$.
Suppose the latter holds.  If the induced subgraph  $[\Delta_i]$
contains an edge, then $[\Delta_i]\cong K_{p+1}$, as the $G$-action
on $p+1$ points is 2-transitive. It follows that $\Gamma=2\cdot
K_{p+1}$ contradicting the fact that $\Gamma$ is connected. Hence
$[\Delta_i]$ does not contain edges of $\Gamma$, and so $\Gamma$ is
a bipartite graph, again a contradiction. Thus $G$ is transitive on
$V\Gamma$.

Let $B_1$ be a block and $u\in B_1$. Then $G_{B_1}\cong
Z_p:Z_{\frac{p-1}{2}}$ and $G_u\cong Z_p:Z_{\frac{p-1}{4}}$. As
$G_u$ has an element of order $p$, $G_u$ is transitive on
$\Gamma(u)$, and hence $\Gamma$ is $G$-arc transitive.

Let $p=5$. Suppose $B_1=\{u,u'\}$.  Since $\Gamma$ is $G$-arc
transitive, it follows that $G_u$ is transitive on $\Gamma(u)$ and
$G_{u'}$ is transitive on $\Gamma(u')$. As  $G_u=G_{u,u'}=G_{u'}
\cong Z_5$ and $\Gamma_3(u)=\{u'\}$, it follows that $\Gamma$ is
$G$-distance transitive. Thus by \cite[p.222, Theorem 7.5.3
(ii)]{BCN}, $\Gamma$ is the icosahedron, so $\Gamma \in
\mathcal{C}(5)$.

Now assume  that $p>5$. As $\Gamma$ is connected and $G$-arc
transitive, $\Gamma\cong Cos(G,H,HgH)$ for the subgroup $H=G_u$ and
some element $g\in G\setminus H$, such that $\langle H,g\rangle =G$
and $g^2\in H$.   Let $a\in H$ and $o(a)=p$. Then $\langle a\rangle$
is a Sylow $p$-subgroup of $G$. Thus $H=\langle a\rangle:\langle
b^2\rangle$ where $N_{G}(\langle a\rangle)=\langle a\rangle:\langle
b\rangle$.

Now we determine the element $g$. Let $u=H$ and $v=Hg$ in $V\Gamma$.
Then $G_u=H$ and $G_{u,v}=\langle b^2\rangle$. Further,
$G_{u,v}^g=(G_u\cap G_v)^g=G_u^g\cap G_v^g=G_v\cap G_u=G_{u,v}$, and
hence $\langle b^2\rangle^g=\langle b^2\rangle$. Thus $g\in
N_G(\langle b^2\rangle)\cong D_{p-1}=\langle b\rangle:\langle
x\rangle$ for some involution $x$. If $g=b^i$ for some $i\geq 1$,
then $\langle H,g\rangle \leq N_X(\langle a\rangle)= \langle
a\rangle:\langle y\rangle$ where $X=PGL(2,p)$ and $y^2=b$,
contradicting the fact that $\langle H,g\rangle =G$. Thus $g=b^ix$
for some $i$, and so $N_G(\langle b^2\rangle)\cong D_{p-1}=\langle
b\rangle:\langle g\rangle$. Thus $\Gamma\cong Cos(G,H,HgH)\in
\mathcal{C}(p)$.

Conversely, assume that $\Gamma \in \mathcal{C}(p)$. If $\Gamma$ is
the icosahedron, then we easily see that $\Gamma$ is a connected
non-bipartite antipodal double cover of $K_6$ and its automorphism
group is $PSL(2,5)\times Z_2$. If $p>5$, then by Remark
\ref{smallval-con-rem}, $\Gamma$ is a connected non-bipartite
antipodal double cover of $K_{p+1}$ and $\Aut \Gamma \cong
PSL(2,p)\times Z_2$.

(b) The claims in part (b) hold for the icosahedron, so assume that
$p>5$   and $p\equiv 1\pmod{4}$,   and let $G=PSL(2,p)$. Let
elements $a_i,b_i,g_i$ and subgroups $H_i$ be chosen as in
Definition \ref{coset-constr} for $i\in \{1,2\}$. Let
$X=PGL(2,p)\cong \Aut G$.

Since all subgroups of $G$ of order $p$ are conjugate there exists
$x\in G$ such that $\langle a_2\rangle^x=\langle a_1\rangle$, so we
may assume that $\langle a_1\rangle=\langle a_2\rangle=M$, say. Let
$Y=N_{X}(M)$. Then $Y=M: \langle y\rangle$ where $o(y)=p-1$, and
$H_1=M:\langle b_1^2\rangle$ and $H_2=M:\langle b_2^2\rangle$ are
equal to the unique subgroup of $Y$ of order $\frac{p(p-1)}{4}$,
that is, $H_1=H_2=M:\langle y^4\rangle=H$, say. Next, since all
subgroups of $Y$ of order $\frac{p-1}{4}$ are conjugate, there exist
$x_1,x_2\in Y$ such that $\langle b_1^2\rangle^{x_1}=\langle
b_2^2\rangle^{x_2}=\langle y^4\rangle$. Since each $x_i$ normalises
$H$ we may assume in addition that $\langle b_1^2\rangle=\langle
b_2^2\rangle=\langle y^4\rangle < \langle y\rangle$. Thus $g_1,g_2$
are non-central involutions in $N_G(\langle y^4\rangle)\cong
D_{p-1}$, an index 2 subgroup of $N_X(\langle y^4\rangle)=\langle
y\rangle:\langle z\rangle\cong D_{2(p-1)}$. The set of non-central
involutions in $N_G(\langle y^4\rangle)$ forms a conjugacy class of
$N_X(\langle y^4\rangle)$ of size $\frac{p-1}{2}$ and consists of
the elements $y^{2i}z$, for $0\leq i<\frac{p-1}{2}$. The group
$\langle y\rangle$ acts transitively on this set of involutions by
conjugation (and normalises $H$). Hence, for some $u\in \langle
y\rangle$, $H^u=H$ and $g_2^u=g_1$. Thus  all graphs in
$\mathcal{C}(p)$ are isomorphic.   Finally, by Lemma
\ref{smallval-p}, these graphs are geodesic transitive of  diameter
3. $\Box$

\section{ Proof of Theorem \ref{gt-primeval3}}

In this section, we will prove Theorem \ref{gt-primeval3} in a
series of lemmas. For all lemmas of this section, we assume that
$\Gamma$ is a connected $2$-geodesic transitive graph of prime
valency $p$ and we denote  $\Aut \Gamma$ by $A$. Note that the
assumption of 2-geodesic transitivity implies that the graph is not
complete. If $\Gamma$ is 2-arc transitive, there is nothing to
prove, so we assume further that this is not the case, that is to
say, we assume that $\Gamma$ has girth 3. The first lemma
determines some intersection parameters.

\begin{lemma}\label{val-p-lemma-1}
Let $(v,u,w)$ be a $2$-geodesic  of $\Gamma$. Then $p\equiv 1 \pmod
4$, $|\Gamma(v)\cap \Gamma(u)|=|\Gamma_2(v)\cap
\Gamma(u)|=\frac{p-1}{2}$ and $|\Gamma(v)\cap \Gamma(w)|$ divides
$\frac{p-1}{2}$. Moreover, $A_v^{\Gamma(v)} \cong
Z_p:Z_{\frac{p-1}{2}}$ is a Frobenius group, and
$A_{v,u}^{\Gamma(v)} \cong Z_{\frac{p-1}{2}}$  is transitive on
$\Gamma(v)\cap \Gamma(u)$.
\end{lemma}
{\bf Proof.} Since $\Gamma$ is $2$-geodesic transitive but not 2-arc
transitive, it follows that  $\Gamma$ is not a cycle. In particular,
$p$ is an odd prime. Let  $|\Gamma(v)\cap \Gamma(u)|=x$ and
$|\Gamma_2(v)\cap \Gamma(u)|=y$. Then $x+y=|\Gamma(u)\setminus
\{v\}|=p-1$. Since   $\girth(\Gamma)=3$, $x\geq 1$. Since $p$ is odd
and the induced subgraph $[\Gamma(v)]$ is an undirected regular
graph with $\frac{px}{2}$ edges,  it follows that $x$ is even. This
together with  $x+y=p-1$ and the fact that $p-1$ is even, implies
that $y$ is also even.

Since $\Gamma$ is arc transitive, $A_v^{\Gamma(v)}$ is transitive on
$\Gamma(v)$.  Since $p$ is a prime, $A_v^{\Gamma(v)}$ acts
primitively on $\Gamma(v)$. By Lemma \ref{val-p-prim-1}, either
$A_v^{\Gamma(v)}$ is 2-transitive, or $A_v^{\Gamma(v)}$ is solvable
and $A_v^{\Gamma(v)}\leq AGL(1,p)$. Since $\Gamma$ is not complete,
it follows that $[\Gamma(v)]$ is not a  complete graph. Also since
$\girth(\Gamma)=3$, $[\Gamma(v)]$ is not an empty graph and so
$A_v^{\Gamma(v)}$ is not 2-transitive. Hence $A_v^{\Gamma(v)}<
AGL(1,p)$. Thus $A_v^{\Gamma(v)}\cong Z_p:Z_m$ is a Frobenius group,
where $m|(p-1)$ and $m<p-1$. Hence $m\leq \frac{p-1}{2}$.

Since $\Gamma$ is vertex transitive, it follows that
$A_u^{\Gamma(u)}\cong Z_p:Z_m$, and hence $A_{u,v}^{\Gamma(u)}\cong
Z_m$ is semiregular on $\Gamma(u)\setminus \{v\}$ with orbits of
size $m$. Since  $\Gamma$ is $2$-geodesic transitive,
$A_{u,v}^{\Gamma(u)}$ is transitive on $\Gamma_2(v)\cap \Gamma(u)$,
and hence  $y=|\Gamma_2(v)\cap \Gamma(u)|=m$, so
$x=p-1-m=m(\frac{p-1}{m}-1)\geq m$, and $x$ is divisible by $m$.

Now  again by arc transitivity, $|\Gamma(u)\cap
\Gamma(w)|=|\Gamma(u)\cap \Gamma(v)|=x$. Since $|\Gamma_2(v)\cap
\Gamma(u)|=m$, it follows that $|\Gamma_2(v)\cap \Gamma(u)\cap
\Gamma(w)|\leq m-1$. Since $\Gamma(w)\cap \Gamma(u)=(\Gamma(w)\cap
\Gamma(u)\cap \Gamma(v))\cup (\Gamma(w)\cap \Gamma(u)\cap
\Gamma_2(v))$,  it follows that
$$x\leq |\Gamma(w)\cap \Gamma(u)\cap \Gamma(v)|+(m-1). \ \ \ \ \ \ \  (*)$$

Let $z=|\Gamma(v)\cap \Gamma(w)|$ and $n=|\Gamma_2(v)|$. Since
$\Gamma$ is $2$-geodesic transitive, $z,n$ are independent of $v,w$
and, counting edges between $\Gamma(v)$ and $\Gamma_2(v)$ we have
$pm=nz$. Now $z\leq |\Gamma(v)|=p$. Suppose first that $z=p$. Then
$m=n$ and $\Gamma(v)=\Gamma(w)$, and so for distinct $w_1,w_2\in
\Gamma_2(v)$, $d_{\Gamma}(w_1,w_2)=2$. Since $\Gamma$ is 2-geodesic
transitive, it follows that $\Gamma(v)=\Gamma(v')$ whenever
$d_{\Gamma}(v,v')=2$. Thus $\diam(\Gamma)=2$, $V\Gamma=\{v\}\cup
\Gamma(v)\cup \Gamma_2(v)$ and $|V\Gamma|=1+p+m$. Let
$\Delta=\{v\}\cup \Gamma_2(v)$. Then for distinct $v_1,v_1'\in
\Delta $, $d_{\Gamma}(v_1,v_1')=2$; for any $v_1''\in
V\Gamma\setminus \Delta$, $v_1,v_1''$ are adjacent. Thus, for any
$v_1\in \Delta$, $\Delta=\{v_1\}\cup \Gamma_2(v_1)$. It follows that
$\Delta$ is a block of imprimitivity for  $A$ of size $m+1$. Hence
$(m+1)|(p+m+1)$, so $(m+1)|p$. Since $m|(p-1)$,  it follows that
$m+1=p$ which contradicts the inequality $m\leq \frac{p-1}{2}$.

Thus $z<p$, and so  $z$ divides $m$, as $pm=nz$.  Since $|\Gamma(w)\cap
\Gamma(u)\cap \Gamma(v)|\leq z$, it follows from  $(*)$ that $x\leq
z+(m-1)\leq 2m-1<2m$.  Since $x$ is divisible by $m$ and $x\geq m$
we have $x=m$. Thus $2m=x+y=p-1$, so $x=y=m=\frac{p-1}{2}$, and
since $x$ is even, $p\equiv 1 \pmod 4$. Also $x=m$ implies that
$A_{v,u}^{\Gamma(v)}$ is transitive on $\Gamma(v)\cap \Gamma(u)$.
Finally, since $nz=pm=p(\frac{p-1}{2})$ and $z<p$, it follows that
$z$ divides $\frac{p-1}{2}$.  $\Box$

\begin{lemma}\label{val-p-lemma-2}
For $v\in V\Gamma$,  the stabiliser $A_v\cong Z_p:Z_{\frac{p-1}{2}}$
is a Frobenius group.
\end{lemma}
{\bf Proof.} Suppose that  $(v,u)$ is an arc of $\Gamma$.  Then by
Lemma \ref{val-p-lemma-1}, $A_v^{\Gamma(v)}\cong
Z_p:Z_{\frac{p-1}{2}}$ is a Frobenius group, and
$A_{v,u}^{\Gamma(v)} \cong Z_{\frac{p-1}{2}}$ is regular on
$\Gamma(v)\cap \Gamma(u)$. Let $K$ be  the kernel of the action of
$A_v$ on $\Gamma(v)$.  Let $u'\in \Gamma(v)\cap \Gamma(u)$ and $x\in
K$. Then $x\in A_{v,u,u'}$. Since $A_{u,v}^{\Gamma(u)} \cong
Z_{\frac{p-1}{2}}$ is semiregular on $\Gamma(u)\setminus \{v\}$, it
follows that $x$ fixes all vertices of $\Gamma(u)$. Since $x$ also
fixes all vertices of $\Gamma(v)$, this argument for each $u\in
\Gamma(v)$ shows that  $x$ fixes all vertices of $\Gamma_2(v)$.
Since $\Gamma$ is connected, $x$ fixes all vertices of $\Gamma$, and
hence $x=1$. Thus $K=1$, so $A_v \cong Z_p:Z_{\frac{p-1}{2}}$ is a
Frobenius group. $\Box$

\begin{lemma}\label{val-p-lemma-3}
Let $(v,u,w)$ be a $2$-geodesic of $\Gamma$. Then $|\Gamma(v)\cap
\Gamma(w)|=\frac{p-1}{2}$, $|\Gamma_2(v)\cap \Gamma(w)\cap
\Gamma(u)|=\frac{p-1}{4}$,  $|\Gamma_2(v)|=p$, and $|\Gamma_2(v)\cap
\Gamma(w)|=\frac{p-1}{2}$.
\end{lemma}
{\bf Proof.}  Let $z=|\Gamma(v)\cap \Gamma(w)|$ and
$n=|\Gamma_2(v)|$. By  Lemma \ref{val-p-lemma-1}, $|\Gamma(u)\cap
\Gamma_2(v)|=\frac{p-1}{2}$ and $z|\frac{p-1}{2}$. Counting the
edges between $\Gamma(v)$ and $\Gamma_2(v)$ gives
$\frac{p-1}{2}p=nz$. By Lemma \ref{val-p-lemma-2},
$A_{v,u}=Z_{\frac{p-1}{2}}$, and by Lemma \ref{val-p-lemma-1},
$A_{v,u}$ is transitive on $\Gamma(v)\cap \Gamma(u)$, so
$[\Gamma(u)]$ is $A_{u}$-arc transitive. Since $p$ is a prime, it
follows by Lemma \ref{arctrans-order-p-1} that $[\Gamma(u)]$ is a
Paley graph $P(p)$. Since $v,w\in \Gamma(u)$ are not adjacent, by
Lemma \ref{val-p-rem1}, $|\Gamma(v)\cap \Gamma(u)\cap
\Gamma(w)|=\frac{p-1}{4}$, hence $z\geq \frac{p-1}{4}+1$. Since
$z|\frac{p-1}{2}$, it follows that $z=\frac{p-1}{2}$. Hence $n=p$.
Thus, $|\Gamma(v)\cap \Gamma(w)|=\frac{p-1}{2}$ and
$|\Gamma_2(v)|=p$.

By Lemma \ref{val-p-lemma-1}, we have $|\Gamma(v)\cap
\Gamma(u)|=\frac{p-1}{2}$. Since $\Gamma$ is arc transitive, it
follows that  $|\Gamma(v_1)\cap \Gamma(v_2)|=\frac{p-1}{2}$ for
every arc $(v_1,v_2)$. Thus, $|\Gamma(u)\cap
\Gamma(w)|=\frac{p-1}{2}$. Since $\Gamma(u)\cap
\Gamma(w)=(\Gamma(v)\cap \Gamma(u)\cap \Gamma(w) )\cup
(\Gamma_2(v)\cap \Gamma(u)\cap \Gamma(w) )$ where  $\Gamma(v)\cap \Gamma(u)\cap \Gamma(w) $ and $\Gamma_2(v)\cap \Gamma(u)\cap \Gamma(w) $ are disjoint, and since $|\Gamma(v)\cap
\Gamma(u)\cap \Gamma(w)|=\frac{p-1}{4}$, it follows that
$|\Gamma_2(v)\cap \Gamma(u)\cap
\Gamma(w)|=\frac{p-1}{2}-\frac{p-1}{4}=\frac{p-1}{4}$.  Since $A_v=
Z_p:Z_{\frac{p-1}{2}}$, it follows that $A_{v,w}=Z_{\frac{p-1}{2}}$
and $A_{v,w}$  is semiregular on $\Gamma_2(v)\setminus \{w\}$ with
orbits of size $\frac{p-1}{2}$. Since $\Gamma_2(v)\cap
\Gamma(w)\subseteq \Gamma(w)\setminus \Gamma(v)$ (of size
$\frac{p-1}{2}$) and since  $|\Gamma_2(v)\cap \Gamma(w)\cap
\Gamma(u)|=\frac{p-1}{4}>0$,  it follows that $|\Gamma_2(v)\cap
\Gamma(w)|=\frac{p-1}{2}$. $\Box$

\begin{lemma}\label{val-p-lemma-4}
Let $v$ be a vertex of $\Gamma$. Then $|\Gamma_3(v)|=1$ and
$\diam(\Gamma)=3$, so $\Gamma$ is antipodal with fibres of size $2$.
Further, $\Gamma$ is geodesic transitive.
\end{lemma}
{\bf Proof.} Suppose that $(v,u,w)$ is a $2$-geodesic of $\Gamma$.
Then by Lemma \ref{val-p-lemma-3}, $|\Gamma(v)\cap
\Gamma(w)|=\frac{p-1}{2}$ and $|\Gamma_2(v)\cap
\Gamma(w)|=\frac{p-1}{2}$. Hence $|\Gamma_3(v)\cap
\Gamma(w)|=p-|\Gamma(v)\cap \Gamma(w)|-|\Gamma_2(v)\cap
\Gamma(w)|=1$. Since $\Gamma$ is $2$-geodesic transitive, it follows
that $|\Gamma_3(v)\cap \Gamma(w_1)|=1$ for all $w_1\in \Gamma_2(v)$.
Thus $\Gamma$ is $3$-geodesic transitive.

Let $\Gamma_3(v)\cap \Gamma(w)=\{v'\}$,  $n=|\Gamma_3(v)|$ and
$i=|\Gamma_2(v)\cap \Gamma(v')|$. Counting edges between
$\Gamma_2(v)$ and  $\Gamma_3(v)$,  we have  $p=ni$.  Since
$[\Gamma(w)]$ is a Paley graph and $u,v'\in \Gamma(w)$ are not
adjacent, it follows from Lemma \ref{val-p-rem1} that
$|\Gamma(u)\cap \Gamma(w)\cap \Gamma(v')|=\frac{p-1}{4}$. Since
$\Gamma(u)\cap \Gamma_2(v)$ contains these $\frac{p-1}{4}$ vertices
as well as $w$, we have $i\geq \frac{p+3}{4}>1$. Thus  $i=p$ and
$n=1$, that is, $|\Gamma_3(v)|=1$. Since $|\Gamma_2(v)\cap
\Gamma(v')|=p$ and $|\Gamma_2(v)|=p$, it follows that
$\Gamma_2(v)=\Gamma(v')$,  so  $\diam(\Gamma)=3$ and $\Gamma$ is
antipodal with fibres of size 2. Therefore $\Gamma$ is geodesic
transitive. $\Box$

\bigskip

We are ready to  prove Theorem \ref{gt-primeval3}.

\bigskip

\noindent {\bf Proof of Theorem \ref{gt-primeval3}.} Let $\Gamma$ be
a connected non-complete graph of prime valency $p$. Suppose first
that $\Gamma$ is 2-geodesic transitive. If $\girth(\Gamma)\geq 4$,
then every 2-arc is a 2-geodesic, so $\Gamma$ is 2-arc transitive.
Now assume that $\girth(\Gamma)=3$. Let $v\in V\Gamma$. Then it
follows from Lemmas \ref{val-p-lemma-1} to \ref{val-p-lemma-4} that
$p\equiv 1 \pmod 4$, $|\Gamma_2(v)|=p$, $|\Gamma_3(v)|=1$ and
$\diam(\Gamma)=3$. Thus, $V\Gamma=\{v\}\cup \Gamma(v)\cup
\Gamma_2(v)\cup \{v'\}$, where $\Gamma_3(v)=\{v'\}$,
$\Gamma(v)=\Gamma_2(v')$ and $\Gamma_2(v)=\Gamma(v')$, and also
$|V\Gamma|=2p+2$. Further, by Lemma \ref{val-p-lemma-4}, $\Gamma$ is
antipodal and geodesic transitive.

Let $\mathcal{B}=\{\Delta_1,\Delta_2, \ldots,\Delta_{p+1}\}$ where
$\Delta_i=\{u_i,u_i'\}$ such that $d_{\Gamma}(u_i,u_i')=3$. Then
each $\Delta_i$ is a block for $A:=\Aut \Gamma$ of size 2 on
$V\Gamma$. Further, for each $j\neq i$, $u_i$ is adjacent to exactly
one vertex of $\Delta_j$, and $u_i'$ is adjacent to the other. The
quotient graph $\Sigma=\Gamma_{\mathcal{B}}$  is therefore a
complete graph $K_{p+1}$ and $\Gamma$ is a cover of $\Sigma$. In
particular, the map $\sigma$  such that $u_i^{\sigma}=u_i'$ and
$u_i'^{\sigma}=u_i$ for all $i$ is an automorphism of $\Gamma$ of
order $2$, and fixes each of the $\Delta_i$ setwise.

We now determine the automorphism group $A$. By Lemma
\ref{val-p-lemma-2}, $A_{v}\cong Z_p{:}Z_{\frac{p-1}{2}}$ is a
Frobenius group, and so $|A|=|A_{v}|.|V\Gamma|=p(p+1)(p-1)$. Let $K$
be the kernel of $A$ acting on ${\mathcal B}$. Then $A$ is an
extension of $K$ by the factor group $A^{{\mathcal B}}$. Since
$\Gamma$ is a cover of $\Sigma$, the kernel $K$ is semiregular on
$V\Gamma$, and hence has order at most 2. Since  the involution
$\sigma$ defined above  lies in $K$, it follows that $K\cong Z_2$.
Thus $|A^{{\mathcal B}}|=|A/K|=\frac{p(p+1)(p-1)}{2}$.

Since $\Gamma$ is arc transitive, the quotient graph
$\Sigma=K_{p+1}$ is $A^{{\mathcal B}}$-arc transitive. Thus,
$A^{{\mathcal B}}$ is 2-transitive on the vertex set ${\mathcal B}$,
and  the point stabiliser $(A^{{\mathcal B}})_{\Delta_1}=
KA_{u_1}/K\cong A_{u_1}\cong Z_p{:}Z_{\frac{p-1}{2}}$ is a Frobenius
group, so $A^{{\mathcal B}}$ is a Zassenhaus group. Since
$|A^{{\mathcal B}}|=\frac{p(p+1)(p-1)}{2}$ and  $A^{{\mathcal B}}$
is not 3-transitive on $\mathcal{B}$, by \cite[Theorem 11.16]{HB-3},
$A^{{\mathcal B}}\cong PSL(2,p)$. Therefore, we have
\[A=K.A^{{\mathcal B}}=Z_2.\PSL(2,p).\]
Suppose that the extension of $Z_2$ by $PSL(2,p)$ is non-split. Then
$A=SL(2,p)$ has only one involution, which lies in the center of
$A$. However, the stabiliser $(A^{{\mathcal B}})_{\Delta_1}\cong
Z_p:Z_{\frac{p-1}{2}}$ is of even order and has trivial center,
which is a contradiction. So the extension $K.A^{{\mathcal B}}$ is
split, and $A\cong Z_2\times PSL(2,p)$. It now follows  from Theorem
\ref{gt-primeval2} (a) that $\Gamma \in \mathcal{C}(p)$.

Conversely, if $\Gamma$ is 2-arc transitive, then it is 2-geodesic
transitive. If $\Gamma\in \mathcal{C}(p)$, then by Theorem
\ref{gt-primeval2} (b), $\Gamma$ is 2-geodesic transitive. $\Box$


\begin{thebibliography}{hhhh}


\bibitem{Baddeley-1}
R. W. Baddeley,  Two-arc transitive graphs and twisted wreath
products. {\it J. Algebraic Combin.} {\bf 2} (1993), 215--237.

\bibitem{Biggs-1}
N. L. Biggs, Algebraic Graph Theory, Cambridge University Press, New
York, (1974).


\bibitem{BCN}
A. E. Brouwer, A. M. Cohen and A. Neumaier, Distance-Regular Graphs,
Springer Verlag, Berlin, Heidelberg, New York, (1989).






\bibitem{DJLP-2}
A. Devillers, W. Jin, C. H. Li and C. E. Praeger, Line graphs and
geodesic transitivity, {\it Ars Math. Contemp.} {\bf 6} (2013),
13--20.



\bibitem{DM-1}
J. D. Dixon and B. Mortimer,  Permutation groups, Springer, New
York, (1996).











\bibitem{GLP}
C. D. Godsil, R. A. Liebler and C. E. Praeger,  Antipodal distance
transitive covers of complete graphs, {\it European J. Combin.} {\bf
19} (1998), 455--478.


\bibitem{GR}
C. D. Godsil and G. F. Royle, Algebraic Graph Theory, Springer, New
York, Berlin, Heidelberg,  (2001).

\bibitem{HB-3}
B. Huppert and N. Blackburn, Finite Groups III, Springer, New York,
Berlin, Heidelberg, (1982).

\bibitem{IP-1}
A. A. Ivanov and C. E. Praeger,  On finite affine 2-arc transitive
graphs. {\it European J. Combin.} {\bf 14} (1993), 421--444.

\bibitem{JLX-val6}
W. Jin, W. J. Liu and S. J. Xu, Finite 2-geodesic transitive graphs of valency 6, submitted.

\bibitem{Li-ci-soluble}
C. H. Li, Finite CI-groups are soluble, {\it Bull. London Math.
Soc.} {\bf 31} (1999), 419--423.




\bibitem{Lorimer-1}
P. Lorimer, Vertex transitive graphs: symmetric graphs of prime
valency, {\it J. Graph Theory} {\bf 8} (1984), 55--68.

\bibitem{Paley-1}
R. E. A. C. Paley,  On orthogonal matrices, {\it J. Math. Phys.}
{\bf 12} (1933), 311--320.

\bibitem{Praeger-4}
C. E. Praeger,  An O'Nan Scott theorem for finite quasiprimitive
permutation groups and an application to 2-arc transitive graphs,
{\it J. London Math. Soc.} {\bf (2) 47} (1993), 227--239.





\bibitem{Taylor-1}
D. E. Taylor,  Two-graphs and doubly transitive groups, {\it J.
Combin. Theory A} {\bf  61} (1992), 113--122.


\bibitem{Tutte-1}
W. T. Tutte,  A family of  cubical  graphs, {\it Proc. Cambridge
Philos. Soc.} {\bf 43} (1947), 459--474.


\bibitem{Tutte-2}
W. T. Tutte, On the symmetry of cubic graphs, {\it Canad. J. Math.}
{\bf 11} (1959), 621--624.

\bibitem{Weiss-1}
R. Weiss,  s-transitive graphs, {\it Colloquia Mathematica
Societatis Janos Bolyai, Algebraic methods in graph theory, szeged
(Hungary)} {\bf 25} (1978), 827--847.

\bibitem{weiss}
R. Weiss,  The non-existence of 8-transitive graphs, {\it
Combinatorica} {\bf 1} (1981), 309--311.






\end{thebibliography}
\end{document}